\documentclass{article} 

\usepackage{amsmath}  
\usepackage{amsthm}
\usepackage[T1]{fontenc}

\theoremstyle{plain}
\newtheorem{thm}{Theorem}
\newtheorem*{mainthm}{Theorem}
\newtheorem{lem}[thm]{Lemma}
\newtheorem{cor}[thm]{Corollary}
\newtheorem{prop}[thm]{Proposition}

\theoremstyle{definition}
\newtheorem{defn}{Definition}

\usepackage{tikz}

\newcommand{\pentRRRRR}{
  \begin{tikzpicture}[baseline=-0.5ex, scale=0.17]
    \draw[red, thick]       (90:1)  -- (18:1);
    \draw[red, thick]      (18:1)  -- (-54:1);
    \draw[red, thick] (-54:1) -- (-126:1);
    \draw[red, thick]    (-126:1) -- (162:1);
    \draw[red, thick]    (162:1) -- (90:1);
  \end{tikzpicture}
}

\newcommand{\pentBRRRR}{
  \begin{tikzpicture}[baseline=-0.5ex, scale=0.17]
    \draw[blue, thick]       (90:1)  -- (18:1);
    \draw[red, thick]      (18:1)  -- (-54:1);
    \draw[red, thick] (-54:1) -- (-126:1);
    \draw[red, thick]    (-126:1) -- (162:1);
    \draw[red, thick]    (162:1) -- (90:1);
  \end{tikzpicture}
}

\newcommand{\pentBRBRR}{
  \begin{tikzpicture}[baseline=-0.5ex, scale=0.17]
    \draw[blue, thick]       (90:1)  -- (18:1);
    \draw[red, thick]      (18:1)  -- (-54:1);
    \draw[blue, thick] (-54:1) -- (-126:1);
    \draw[red, thick]    (-126:1) -- (162:1);
    \draw[red, thick]    (162:1) -- (90:1);
  \end{tikzpicture}
}

\newcommand{\pentBBRRR}{
  \begin{tikzpicture}[baseline=-0.5ex, scale=0.17]
    \draw[blue, thick]       (90:1)  -- (18:1);
    \draw[blue, thick]      (18:1)  -- (-54:1);
    \draw[red, thick] (-54:1) -- (-126:1);
    \draw[red, thick]    (-126:1) -- (162:1);
    \draw[red, thick]    (162:1) -- (90:1);
  \end{tikzpicture}
}

\title{Ramsey Theory for Edge Ideals}
\author{ Alexander Engström\\ alexander.engstrom@fra.se }
\date\today

\begin{document}
\maketitle

\begin{abstract}
We initiate the study of Ramsey theory for edge ideals. The edge ideal $I_G$ is generated by the monomials $x_ux_v$ for all edges $uv$ of $G.$ A central invariant measuring structure in algebraic geometry and commutative algebra are the betti number $\beta_{i,j}(S/I)$ for some ideal $I$ in a ring $S.$ We prove the Ramsey Theorem for Edge Ideals: For any betti number indices $i,j,$ and any number of colors $m,$ there is a graph $F$ such that for every $m$-coloring of the edges of $F,$ there is a color $l$ such that $ \beta_{i,j}(S/I_{F_l}) >0, $ where $F_l$ is the subgraph of $F$ whose edges are in color $l$. Moreover, the number of vertices of $F$ is at most $jm^2.$
\end{abstract}

\section{Introduction}

Ramsey theory asserts that any sufficiently large graph, regardless of how its edges are partitioned or colored, must contain structured monochromatic substructures. A classical illustration of this principle is that any 2-coloring of the edges of the complete graph $K_6$ necessarily yields a monochromatic triangle. During the one hundred years of Ramsey theory, its deepest results are among the hallmarks of combinatorics in its interaction with other parts of mathematics, for example the existence of arbitrary long arithmetic progressions in the primes.

In this paper we initiate the exploration of synergies between Ramsey theory and commutative algebra, in particular edge ideals. For a graph $G$ consider the polynomial ring $S$ whose variables are indexed by the vertices of $G.$ The edge ideal $I_G$ is generated by the monomials $x_ux_v$ for all edges $uv$ of $G.$ For an ideal $I,$ the central invariants in commutative algebra encoding the algebraic geometric information are the betti numbers $\beta_{i,j}(S/I),$ which are indexed by a pair $i,j$ of positive integers. The existence of geometric structure is certified by $\beta_{i,j}(S/I)>0.$

\begin{mainthm}[The Ramsey Theorem for Edge Ideals]
For any betti number indices $i,j,$ and number of colors $m,$ there is a graph $F$ such that for every $m$-coloring of the edges of $F,$ there is a color $l$ such that
\[ \beta_{i,j}(S/I_{F_l}) >0, \]
where $F_l$ is the subgraph of $F$ whose edges are in color $l$. Moreover, the number of vertices of $F$ is at most $jm^2.$
\end{mainthm}

The exact formulation is in Theorem~\ref{thm:main} and the upper bound in Corollary~\ref{cor:upper}, which are both in Section 4. In that section we also provide the foundational definitions for Ramsey theory for edge ideals, prove some further results and provide some examples. Before that, in Section 2, we explain some useful tools from graph theory, topological combinatorics, and commutative algebra, as well as include plenty of references to papers and textbooks aimed at the interaction between methods within those areas. In Section 3 we give a brief introduction to the typical questions asked in Ramsey theory to provide context, but that part is not necessary for Section 4 and the main theorem of the paper, it is only informative.

\section{Basic notions}

\subsection{Graph theory}
All graphs are finite and undirected. The vertex set of a graph $G$ is denoted by $V(G)$ and the edge set by $E(G)$. A \emph{subgraph} $H \subseteq G$ is a selection of vertices of $G$ and edges among them from $G.$ A subgraph is \emph{induced} if all possible edges among the selected vertices are included. The induced subgraph on $W \subseteq V(G)$ is denoted $G[W].$ Special classes of graphs as cycles $C_n$, paths $P_n$ and complete graphs $K_n$ are indexed by their number of vertices $n.$ The disjoint union of two graphs $G$ and $H$ is denoted $G \sqcup H.$ A subset of vertices $W\subset V(G)$ is \emph{independent} in $G$ if there are no edges in $G[W]$. For further basic graph theory, see Diestel~\cite{diestel2025}.

\subsection{Topological combinatorics}
A \emph{simplicial complex} $\Delta$ on a vertex set $V$ is a set of subsets of $V$ such that if $\tau \in \Delta$ and $\sigma \subseteq \tau,$ then $\sigma \in \Delta.$ The vertex set of all simplicial complexes in this paper are finite. The \emph{reduced homology} of a simplicial complex $\Delta$ in dimension $d$ over a field $\mathbf{k}$ is denoted 
$\widetilde{H}_d(\Delta ; \mathbf{k}).$ Whenever the specific field $\mathbf{k}$ is not relevant, it is dropped. The \emph{join} of two simplicial complexes $\Delta$ and $\Gamma$ on disjoint vertex sets is $\Delta \ast \Gamma = \{\delta \cup \gamma \mid \delta \in \Delta \textrm{ and } \gamma \in \Gamma \}.$
For further basic algebraic topology, see Matou\v{s}ek~\cite{matousek2008}, Munkres~\cite{munkres1984} or Spanier~\cite{spanier1966}.
\begin{defn}
The \emph{independence complex} of a graph $G,$ denoted by $\operatorname{Ind}(G),$ is the simplicial complex of the independent sets of $G.$
\end{defn}
For further results on independence complexes of graphs, see for example Bayer et al~\cite{bayer2026}, Chudnovsky et al~\cite{chudnovsky2020}, Ehrenborg and Hetyei~\cite{ehrenborg2006}, Engstr\"om~\cite{engstrom2008, engstrom2009}, Jonsson~\cite{jonsson2008, jonsson2010} or Kim~\cite{kim2022, kim2026}.

\subsection{Commutative algebra}
We restrict our presentation of commutative algebra strictly to the tools required for our applications rather than developing the theory in full generality.
\begin{defn}\label{defn:edgeideal}
For a graph $G,$ the \emph{edge ideal} $I(G)$ in the ring $S=\mathbf{k}[x_v \mid v \in V(G) ]$ is generated by the monomials $x_ux_v$ for all edges $uv$ of $G.$
\end{defn}
For an arbitrary ideal $I$ in a ring $S,$ one can encode the successive modules of syzygies—the linear dependencies among generators, the relations among those relations, and so forth. The collection of that data is a \emph{free resolution} and its \emph{length} is the number of times you need relations among relations among relations and so on. For a polynomial ring on $n$ variables, Hilbert's Syzygy Theorem from 1890 states that the length is finite and at most $n$. In a well defined sense, minimal resolutions are unique, see Eisenbud \cite{eisenbud1995}. The dimensions of the pieces of the minimal resolution are encoded in the \emph{betti numbers} $\beta_{i,j}(S/I).$ Those numbers are intrinsically linked to geometrical phenomena of the scheme cut out by $I$ and were introduced to capture intuitive reasoning about intersections between curves in a stringent manner, see Eisenbud~\cite{eisenbud2005}. For edge ideals there is a more direct way to compute their betti numbers, and we take that as its definition. By tensoring with the koszul resolution, or applying Hochster's formula, as explained by Miller and Sturmfels~\cite{miller2005}, we have that:
\begin{defn}\label{defn:betti} For any graph $G,$ set
\[
\beta_{i,j}(S/I(G)) = \sum_{\substack{W \subseteq V(G) \\ |W| = j}} \dim_{\mathbf{k}} \widetilde{H}_{j - i - 1}(\operatorname{Ind}(G[W]); \mathbf{k})
\]
for any positive integers $i$ and $j$ such that $j/2 \leq i \leq j-1.$ 
\end{defn}
Two important concepts of resolutions and betti numbers are  \emph{linearity} and \emph{regularity}, which we will however not define. In a foundational paper on edge ideals, Fr\"oberg~\cite{froberg1990} proved that the resolution is linear exactly when the graph is the complement of a chordal graph. Engstr\"om and Orlich~\cite{engstrom2023} employed probabilistic techniques common in Ramsey theory to show that the regularity of edge ideals is determined for almost all graphs.

There is a vast literature on edge ideals and their relation to combinatorics, see for example Dochtermann and Engstr\"om~\cite{dochtermann2009}, Dochtermann and Newman~\cite{dochtermann2023}, Engstr\"om and Stamps~\cite{engstrom2013}, Ha and Van Tuyl~\cite{ha2008}, Herzog and Hibi~\cite{herzog1999}, Van Tuyl and Villarreal~\cite{vantuyl2008} and the textbook by Herzog and Hibi~\cite{herzog2011}.

\section{Ramsey theory context}
In this section we will briefly touch upon the world of Ramsey theory. None of the content of this section is necessary for the remainder of this paper. It's purely here to contextualize and to explain the choice of notation, definitions, and theorems stated in the next section.

The \emph{classical (or ordinary) Ramsey number} of two graphs $G$ and $H$, denoted by $R(G, H)$, is the minimum positive integer $N$ such that every 2-coloring of the edges of the complete graph $K_N$ contains either a red copy of $G$ or a blue copy of $H$.

The existence of $R(G,H)$ for finite graphs was proved by Ramsey in 1928 \cite{ramsey1930}. For a basic introduction to Ramsey theory, see for example Chapter 9 in Diestel \cite{diestel2025}. There are also Ramsey numbers for $r\geq 2$ colors. If one try to avoid $K_r$ in each color, then the multicolor Ramsey number is at most $e^{-\delta k}r^{rk}$ for some constant $\delta>0$ depending on $k,$ and sufficiently large $k,$ as proved by Balister et al \cite{balister2026}.

The \emph{induced Ramsey number} of two graphs $G$ and $H$, denoted by $IR(G, H)$, is the minimum positive integer $N$ for which there exists a graph $F$ on $N$ vertices such that every 2-coloring of the edges of $F$ contains either a red induced copy of $G$ or a blue induced copy of $H$. In Ramsey theory it is then common to state that $F$ \emph{arrows} $G$ and $H,$
\[ F \rightarrow (G,H), \]
a notation we will lean on later.

The existence of $IR(G,H)$ for all finite graphs $G$ and $H$ was established independently by Deuber~\cite{deuber1975}, Erd\H{o}s, Hajnal, and P\'osa~\cite{erdos1975}, and R\"odl~\cite{rodl1973}.
Conlon, Fox, and Sudakov~\cite{conlon2010} established near-exponential upper bounds for induced Ramsey numbers, significantly improving upon previous tower-type estimates.
The long-standing Kohayakawa--Pr\"omel--R\"odl conjecture~\cite{kohayakawa1998} concerning polynomial bounds for off-diagonal cases was recently resolved in a preprint on the arXiv by Andrades, Campos, Filipe, and Kats~\cite{andrades2026} as $IR(G,H) \leq |V(H) | ^ { 60 [V(G) |^2 +3 }.$

The \emph{weak induced Ramsey number} of two graphs $G$ and $H$, denoted by $IR_w(G, H)$, is the minimum positive integer $N$ for which there exists a graph $F$ on $N$ vertices such that for every 2-coloring of the edges of $F$ (partitioning the edges into a red subgraph $F_r$ and a blue subgraph $F_b$), there is either an induced copy of $G$ in $F_r$ or an induced copy of $H$ in $F_b$.

The weak induced Ramsey number $IR_w(G,H)$ was introduced by Dudek, Grytczuk, and Ruci\'nski~\cite{dudek2023} and further investigated in~\cite{debski2023}.

Another direction of generalization is Ramsey theory of properties of graphs or sets of graphs. Given two finite sets of graphs $\mathcal{G}$ and $\mathcal{H},$ the Ramsey number $R(\mathcal{G},\mathcal{H}) \leq N$ whenever a two coloring of $K_N$ contains a red $G\in  \mathcal{G}$ or a blue $H \in \mathcal{H}.$ Using the terminology of weak induced Ramsey theory for sets of graphs, the Ramsey numbers of Betti numbers are nothing else than defined by the sets of graphs with certain reduced homology as stated in Definition~\ref{defn:betti}.

\section{Ramsey theory for edge ideals}

\begin{defn}[Ramsey arrow notation for betti numbers]
Let $i_1,\ldots,i_m$ and $j_1,\ldots, j_m$ be positive integers satisfying $j_l/2 \leq i_l \leq j_l-1;$  $\mathbf{k}$ a field and $F$ a graph. If there for every $m$-edge coloring of $F$ (partitioning it into the graphs $F_1, F_2, \ldots, F_m$) is an $l$ such that $\beta_{i_l,j_l}^{\mathbf{k}}(S/I(F_l))>0,$ then 
\[ F \xrightarrow{\beta} (i_1,j_1;i_2,j_2;\ldots;i_m,j_m; \mathbf{k}). \]
\end{defn}

The existence of an $F$ follows the weak induced Ramsey theorem mentioned in the previous section. In Theorem~\ref{thm:main} we will prove a stronger theorem based on elementary considerations.

\begin{thm}\label{thm:arrowsExist}
For any positive integers $i_1,\ldots,i_m$ and $j_1,\ldots, j_m$ satisfying $j_l/2 \leq i_l \leq j_l-1$ and a field $\mathbf{k},$ there is a graph $F$ such that
\[ F \xrightarrow{\beta} (i_1,j_1;i_2,j_2;\ldots;i_m,j_m; \mathbf{k}). \]
\end{thm}
\begin{proof}
For an $m$-edge coloring of $F,$ we need to ensure that for every $m$-edge coloring, $\beta_{i_l,j_l}^{\mathbf{k}}(S/I(F_l))>0$ for some $l.$ By Definition~\ref{defn:betti}, this equates to that
$\dim_{\mathbf{k}} \widetilde{H}_{j_l - i_l - 1}(\operatorname{Ind}(F_l[W]); \mathbf{k})>0$ for some $W \subseteq V(F)$ on $j$ vertices. For each $l,$ pick a graph $G_l$ on $j_l$ vertices such that
$\dim_{\mathbf{k}} \widetilde{H}_{j_l - i_l - 1}(\operatorname{Ind}(G_l); \mathbf{k})>0.$ By the weak induced Ramsey theorem, there is a graph $F$ such that any $m$-edge coloring of it will have an induced monochromatic $G_l$ for some color $l.$
\end{proof}

\begin{defn}[Ramsey number notation for betti numbers] Let $i_1,\ldots,i_m$ and $j_1,\ldots, j_m$ be positive integers satisfying $j_l/2 \leq i_l \leq j_l-1$ and $\mathbf{k}$ a field. The Ramsey number
\[ R_\beta (i_1,j_1;i_2,j_2;\ldots;i_m,j_m; \mathbf{k}) \]
is the smallest positive integer such that there is a graph $F$ on that number of vertices satisfying
\[ F \xrightarrow{\beta} (i_1,j_1;i_2,j_2;\ldots;i_m,j_m; \mathbf{k}). \]
\end{defn}

In ordinary graph theory the hardest cases are when the number of colors grow while the other parameters are fixed, it's called the diagonal case.

\begin{defn}[Diagonal Ramsey number notation for betti numbers] Let $i,j,m$ be positive integers satisfying $j/2 \leq i \leq j-1$ and $\mathbf{k}$ a field. Set
\[ R_\beta^m (i,j; \mathbf{k}) = R_\beta (i_1,j_1;i_2,j_2;\ldots;i_m,j_m; \mathbf{k}) \]
where $i_l = i$ and $j_l = j.$
\end{defn}

Betti numbers do sometimes depend on the choice of the field $\mathbf{k},$ but it is unclear if there are non-trivial examples when Ramsey theory does. Any arguments in this paper will be independent of $\mathbf{k}$ and the notation for how to drop the field has to be agreed upon.

\begin{defn}[Field independent notations] Let $i_1,\ldots,i_m$ and $j_1,\ldots, j_m$ be positive integers satisfying $j_l/2 \leq i_l \leq j_l-1.$ If
$ F \xrightarrow{\beta} (i_1,j_1;i_2,j_2;\ldots;i_m,j_m; \mathbf{k}) $ for any field $\mathbf{k},$ then $ F \xrightarrow{\beta} (i_1,j_1;i_2,j_2;\ldots;i_m,j_m).$
Set $R_\beta (i_1,j_1;i_2,j_2;\ldots;i_m,j_m) = \max_{\mathbf{k}} R_\beta (i_1,j_1;i_2,j_2;\ldots;i_m,j_m; \mathbf{k}).$ Analogously for $R_\beta^m (i,j).$
\end{defn}

That the $\max_{\mathbf{k}}$ in the definition is well defined follows from Theorem~\ref{thm:main} or by selecting the $G_l$ in the proof of Theorem~\ref{thm:arrowsExist} as disjoint unions of complete graphs.
We need two elementary lemmas from algebraic topology, the first one follows right off from definitions and the other one from the K\"unneth formula for joins.
\begin{lem}\label{lem:compH0} 
The number of connected components of the complement of a graph $G$ is $\dim_K  \widetilde{H}_0(  \operatorname{Ind}(G) ; \mathbf{k} )+1.$
\end{lem}
\begin{lem}\label{lem:join} 
If $ \dim_{\mathbf{k} }  \widetilde{H}_{d_s}(\Delta_s; {\mathbf{k} }) >0$ for $s=1,2,\ldots,t$ then
\[ \dim_{\mathbf{k} }  \widetilde{H}_{d_1 + d_2 + \cdots d_t + t-1}(\Delta_1 \ast \Delta_2 \ast \cdots \ast \Delta_t; {\mathbf{k} }) >0.\]
\end{lem}
Combining these lemmas yields the technical workhorse required for our main theorem:
\begin{lem}\label{lem:horse} If $G=H_1\sqcup H_2 \sqcup \cdots \sqcup H_t,$ and for each graph $H_s$ 
\begin{itemize}
\item[-] there are at least two vertices in $H_s$ and
\item[-] one of the vertices of $H_s$ is adjacent to all the other vertices,
\end{itemize}
then $\dim_{\mathbf{k} }  \widetilde{H}_{t-1}(  \operatorname{Ind}(G) ; {\mathbf{k} } )>0.$
\end{lem}
\begin{proof}
By Lemma~\ref{lem:compH0} we have that
$\dim_{\mathbf{k} }  \widetilde{H}_{0}(  \operatorname{Ind}(H_s) ; {\mathbf{k} } )>0$ for each $s,$ as the vertex in $H_s$ adjacent to all others is isolated in the complement of $H_s$ from the other vertices. The independence complex of the disjoint union of graphs is the join of their independence complexes. By Lemma~\ref{lem:join} we are done.
\end{proof}
This is our main theorem.
\begin{thm}\label{thm:main}
Let $i_1,\ldots,i_m$ and $j_1,\ldots, j_m$ be positive integers satisfying $j_l/2 \leq i_l \leq j_l-1.$ If
$N=2+\sum_{l=1}^m \left( \left\lceil \frac{j_l}{j_l-i_l} \right\rceil -2 \right)$
and $M=1+\sum_{l=1}^m (j_l-i_l-1)$
then
\[ \bigsqcup_M K_N \xrightarrow{\beta} (i_1,j_1;i_2,j_2;\ldots;i_m,j_m ). \]
\end{thm}
\begin{proof}
Fix an $m$-edge coloring of $\bigsqcup_M K_N$ and a designated vertex of each $K_N.$ Now consider the coloring of each $K_N$ one by one.

There are $N-1=1+\sum_{l=1}^m \left( \left\lceil \frac{j_l}{j_l-i_l} \right\rceil -2 \right)$ edges from the designated vertex to the rest of $K_N.$ By the pigeonhole principle, there has to be a color $l$ with at least 
\[ 1 +  \left( \left\lceil \frac{j_l}{j_l-i_l} \right\rceil -2 \right) =   \left\lceil \frac{j_l}{j_l-i_l} \right\rceil -1 \]
edges from the designated vertex in color $l.$ In that case, we say that the color $l$ \emph{dominates} that $K_N.$

Every $K_N$ is dominated by a color, and there is a disjoint union of $M=1+\sum_{l=1}^m (j_l-i_l-1)$ of them in $\bigsqcup_M K_N.$ By the pigeonhole principle, there is a color $l$ that dominates at least
\[1 + (j_l-i_l-1) = j_l-i_l \] of them.

Fix the two non-negative integers $n_0$ and $n_1$ such that $n_0+n_1= j_l-i_l$ and
\[n_0  \left\lceil \frac{j_l}{j_l-i_l} \right\rceil + n_1 \left\lfloor \frac{j_l}{j_l-i_l} \right\rfloor = j_l.
\]
For $n_0$ of the $K_N$ dominated by the color $l$, select subsets $W_1,W_2, \ldots W_{n_0}$ of cardinality $\left\lceil \frac{j_l}{j_l-i_l} \right\rceil$ of their vertex sets such that
\begin{itemize}
\item[-] the designated vertex is in each $W_s$ and
\item[-] all the edges from the designated vertex within each $W_s$ is in color $l.$
\end{itemize}
Then select $W_{n_0+1},W_{n_0+2}, \ldots W_{j_l-i_l}$ analogously but each of cardinality $\left\lfloor \frac{j_l}{j_l-i_l} \right\rfloor.$ 

Let $F_l$ be the subgraph of $F$ defined by the $l$-colored edges. The cardinality of $W_1 \cup W_2 \cup \ldots \cup W_{j_l-i_l}$ is $j_l$ and by Lemma~\ref{lem:horse}
\[
\dim_{\mathbf{k} }  \widetilde{H}_{j_l-i_l-1}(  \operatorname{Ind}(F_l[W_1 \cup W_2 \cup \ldots \cup W_{j_l-i_l}]) ; {\mathbf{k} } )>0.
 \]
By Definition~\ref{defn:betti} we have that $\beta_{i_l,j_i}(S/I(F_l)) >0$ and 
\[\bigsqcup_M K_N \xrightarrow{\beta} (i_1,j_1;i_2,j_2;\ldots;i_m,j_m ).\]
\end{proof}

\begin{cor}\label{cor:ramseybetti}
For positive integers $i_1,\ldots,i_m$ and $j_1,\ldots, j_m$ satisfying $j_l/2 \leq i_l \leq j_l-1,$ 
\[R_\beta(i_1,j_1;i_2,j_2;\ldots;i_m,j_m) \leq \left( 1+\sum_{l=1}^m (j_l-i_l-1) \right) \left(  2+\sum_{l=1}^m \left( \left\lceil \frac{j_l}{j_l-i_l} \right\rceil -2 \right)  \right). \]
\end{cor}

\begin{cor}\label{cor:multi} For any positive integers $i,j,m$ satisfying $j/2 \leq i \leq j-1,$ 
\[R_\beta^m(i,j) \leq \left( 1+ m (j-i-1) \right) \left(  2(1-m) + m  \left\lceil \frac{j}{j-i} \right\rceil  \right).
\]
\end{cor}

We require the following elementary result:
\begin{prop} If $i,j,m$ are positive integers satisfying $j/2 \leq i \leq j-1$ and $m\geq 2,$ then
\[ 
\left( 1+ m (j-i-1) \right) \left(  2(1-m) + m  \left\lceil \frac{j}{j-i} \right\rceil  \right)   \leq j m^2 \]
\end{prop}
\begin{proof}
\[ 
\begin{array}{rcl}
\left( 1+ m (j-i-1) \right) \left(  2(1-m) + m  \left\lceil \frac{j}{j-i} \right\rceil  \right) &\leq & m (j-i)   \left(  2 -2m + m \left(  \frac{j}{j-i} +1 \right)  \right) \\
& = & m (j-i)   \left(  2 -m  +m  \frac{j}{j-i}  \right) \\
& \leq & m (j-i)  m  \frac{j}{j-i}  \\
& = & jm^2   \\
\end{array}
\]
\end{proof}

\begin{cor}\label{cor:upper} If $i,j,m$ are positive integers satisfying $j/2 \leq i \leq j-1$ and $m\geq 2,$ then $ R_\beta^m(i,j) \leq jm^2. $
\end{cor}

For certain combination of the indices $i$ and $j$ there are linear bounds in $m.$ In commutative algebra the betti numbers are tabulated in a slightly particular way, and the corollaries are named from that.

\begin{cor}[linear strand] $R_\beta^m(i,i+1)  \leq (i-1)m+2.$
\end{cor}
\begin{proof}
$
R_\beta^m(i,i+1)  \leq  \left( 1+ m (i+1-i-1) \right) \left(  2(1-m) + m  \left\lceil \frac{i+1}{i+1-i} \right\rceil  \right) =  (i-1)m+2.
$
\end{proof}
\begin{cor}[diagonal] $R_\beta^m(i,2i)  \leq 2(i-1)m+2.$
\end{cor}
\begin{proof}
$
R_\beta^m(i,j)  \leq  \left( 1+ m (2i-i-1) \right) \left(  2(1-m) + m  \left\lceil \frac{2i}{2i-i} \right\rceil  \right) = 2(i-1)m+2. \\
$
\end{proof}

In Ramsey theory it is common that bounds for small Ramsey numbers lift to higher ones, but usually involving large factors. This is much more moderate.

\begin{thm}\label{thm:lift}
Let $r_1, r_2, \ldots, r_m$ be positive integers and $r=r_1+r_2+\cdots + r_m - m + 1$. If
\[ F \xrightarrow{\beta} (i_1,j_1;i_2,j_2;\ldots;i_m,j_m; \mathbf{k} ) \]
then
\[ \bigsqcup_r F  \xrightarrow{\beta} (r_1i_1,r_1j_1;r_2i_2,r_2j_2;\ldots;r_mi_m,r_mj_m; \mathbf{k})   \]
\end{thm}
\begin{proof} Fix an $m$-edge coloring of $\bigsqcup_r F.$ For each copy of $F$ it is assumed that $F \xrightarrow{\beta} (i_1,j_1;i_2,j_2;\ldots;i_m,j_m; \mathbf{k} )$
and there is a color $l$ certifying that by $\beta_{i_l,j_l}^{\mathbf{k}}(S/I(F_l))>0.$ To be more precise, there is a $W\subseteq F$ of cardinality $j_l$ such that
\[
\dim_{\mathbf{k}} \widetilde{H}_{j_l - i_l - 1}(\operatorname{Ind}(F_l[W]); \mathbf{k}) >0.
\]
In this situation, we declare that the color $l$ dominates that particular copy of $F.$

Every copy of $F$ is dominated by some color. By the pigeonhole principle, there is a color $l$ dominating at least $r_l$ copies of $F.$ Denote the sets certifying the color $l$ dominance by
$W_1, W_2, \ldots, W_{r_l}.$ By definition
\[
\dim_{\mathbf{k}} \widetilde{H}_{j_l - i_l - 1}(\operatorname{Ind}(F_l[W_s]); \mathbf{k}) >0
\]
for each $s=1,2,\ldots, r_l.$ The independence complex of a disjoint union of graphs is the join of their independence complexes. By Lemma~\ref{lem:join},
\[
\dim_{\mathbf{k}} \widetilde{H}_{r_lj_l - r_li_l - 1}(\operatorname{Ind}(F_l[W_1 \sqcup W_2 \sqcup \cdots \sqcup W_{r_l}   ]); \mathbf{k}) >0.
\]
This shows that
\[
\beta_{r_li_l,r_lj_l}(S/I(F_l))>0
\]
and
\[ \bigsqcup_r F  \xrightarrow{\beta} (r_1i_1,r_1j_1;r_2i_2,r_2j_2;\ldots;r_mi_m,r_mj_m; \mathbf{k}).   \]
\end{proof}

\begin{cor}\label{cor:lift} Let $r=r_1+\cdots + r_m -m+1.$
\[R_\beta (r_1i_1,r_1j_1 ;\ldots;r_mi_m,r_mj_m; \mathbf{k})
\leq  r R_\beta (i_1,j_1 ;\ldots;i_m,j_m; \mathbf{k}).
\]
\end{cor}

\begin{cor}
\[R_\beta^m (ri,rj ; \mathbf{k})
\leq  (m(r-1)+1) R_\beta (i,j ; \mathbf{k}).
\]
\end{cor}

\subsection{Explicit calculation example}

It is usually possible to improve on Theorem~\ref{thm:main} by explicit calculations. It states that
\[ \bigsqcup_3 K_4 \xrightarrow{\beta} (3,5;3,5)  \]
as
$2+2 \left( \left\lceil \frac{5}{5-3} \right\rceil -2 \right)=4$ and $1+ 2(5-3-1)=3.$ While Theorem~\ref{thm:main} provides a general upper bound, it can often be refined via explicit constructions. Although a systematic pruning of vertices from $\bigsqcup_3 K_4$ could theoretically sharpen the estimate, the resulting proofs would become cumbersome. Instead, direct explicit analysis yields a tighter bound:
\begin{prop}\label{prop:ex} $C_5 \sqcup K_2 \xrightarrow{\beta} (3,5;3,5)$ and  $R_\beta(3,5;3,5)\leq 7.$
\end{prop}
\begin{proof}
Assume the contradiction, that there is a way to color the edges of $C_5 \sqcup K_2$ in red and blue without an induced monochromatic subgraph on five vertices with homology in the first dimension.

We may assume that the edge coloring have more red than blue edges on the five cycle. There are four cases:
\begin{itemize}
\item[(\pentRRRRR)] The induced red subgraph is $C_5$ and $\dim \widetilde{H}_1(\operatorname{Ind}(C_5)) = 1,$
\item[(\pentBRRRR)] The induced red subgraph is $P_5$ and $\dim \widetilde{H}_1(\operatorname{Ind}(P_5)) = 1,$ 
\item[(\pentBRBRR)] The induced red subgraph is $K_2 \sqcup P_3$ and $\dim \widetilde{H}_1(\operatorname{Ind}(K_2 \sqcup P_3)) = 1,$
\item[(\pentBBRRR)] Both the induced red and blue subgraphs on five vertices have isolated vertices and thus no homology in its independence complexes.
\end{itemize}
The only possible coloring is \pentBBRRR. Note that it includes both a blue and a red induced $P_3.$ In the coloring of $C_5 \sqcup K_2$ we have so far not considered the color of the edge $K_2,$ but we know that there has to be a $P_3$ in the $C_5$ of the same color. We have an induced monochromatic $K_2 \sqcup P_3.$ But $\dim \widetilde{H}_1(\operatorname{Ind}(K_2 \sqcup P_3)) = 1,$ which is a contradiction.
\end{proof}

Returning to Theorem~\ref{thm:main}, as condensed in Corollary~\ref{cor:ramseybetti}, we had that 
\[
\begin{array}{rcl}
R_\beta(3r,5r;3r,5r) & \leq & \left( 1+ 2(5r-3r-1) \right) \left(  2+ 2 \left( \left\lceil \frac{5r}{5r-3r} \right\rceil -2 \right)  \right) \\
& = & 16r-4
\end{array}
\]
for any positive integer $r.$ The result of Proposition~\ref{prop:ex} together with Theorem~\ref{thm:lift} as condensed in Corollary~\ref{cor:lift} proves the slightly stronger result
\[R_\beta (3r,5r;3r,5r)  \leq  (2r-1) R_\beta (3,5;3,5) = 7(2r-1) = 14r-7. \]

\subsubsection*{AI use} Only to edit language, not mathematics.

\end{document}